\documentclass{amsart}

\usepackage{hyperref}
\usepackage{amsthm,mathrsfs,amssymb,amsmath} 
\usepackage{pstricks,pst-node,pst-tree}

\begin{document}

\title{$q$-analog of tableau containment}
\author{Jang Soo Kim}
\thanks{The author was supported by the SRC program of Korea Science and
  Engineering Foundation (KOSEF) grant funded by the Korea government (MEST)
  (No. R11-2007-035-01002-0).}
\email{kimjs@math.umn.edu}

\begin{abstract}
  We prove a $q$-analog of the following result due to McKay, Morse and Wilf:
  the probability that a random standard Young tableau of size $n$ contains a
  fixed standard Young tableau of shape $\lambda\vdash k$ tends to
  $f^{\lambda}/k!$ in the large $n$ limit, where $f^{\lambda}$ is the number of
  standard Young tableaux of shape $\lambda$. We also consider the
  probability that a random pair $(P,Q)$ of standard Young tableaux of the same
  shape contains a fixed pair $(A,B)$ of standard Young tableaux. 
\end{abstract}

\maketitle

\newtheorem{thm}{Theorem}[section]
\newtheorem{lem}[thm]{Lemma}
\newtheorem{prop}[thm]{Proposition}
\newtheorem{cor}[thm]{Corollary}
\theoremstyle{definition}
\newtheorem{defn}{Definition}[section]
\newtheorem{conj}{Conjecture}[section]
\newtheorem{problem}{Problem}[section]

\newcommand{\seq}[3]{#1_1 #3 #1_2 #3 \def\temp{#3} \def\comma{,}
  \if\comma\temp \ldots \else \cdots \fi #3 #1_{#2} }
\newcommand{\mat}[2]{\left(\begin{array}{#1}#2\end{array}\right)}

\newcommand{\sh}{\operatorname{sh}}
\newcommand{\row}{\operatorname{row}}
\newcommand{\col}{\operatorname{col}}
\newcommand{\inv}{\operatorname{inv}}
\newcommand{\maj}{\operatorname{maj}}
\newcommand{\imaj}{\operatorname{imaj}}

\def\qbinom#1#2{\genfrac{[}{]}{0pt}{}{#1}{#2}_q}
\def\pbinom#1#2{\genfrac{[}{]}{0pt}{}{#1}{#2}_p}
\def\biletter#1#2{\genfrac{}{}{0pt}{}{#1}{#2}}
\def\xfactorial#1#2{\left[ #1 \right]_{#2}!}
\def\qfactorial#1{\xfactorial{#1}{q}}
\def\pfactorial#1{\xfactorial{#1}{p}}
\def\limit#1{\lim_{#1\rightarrow \infty}}
\def\pp{\bar p}
\def\qq{\bar q}
\def\rr{\bar r}
\def\sf{\mathsf{sf}}
\def\ss{{\bf s}}
\def\ov{\overline}

\def\vr{{\mathbf r}}
\def\vc{{\mathbf c}}

\def\t{\mathcal{T}}
\def\N{\mathbb{N}}
\def\S{\mathfrak{S}}
\def\I{\mathfrak{I}}

\def\sp{\S_n(\sigma,\tau)}
\def\spl{\S_\ell(\sigma,\tau)}

\def\ul#1{^{\leq {#1}}}
\def\ug#1{^{> {#1}}}
\def\dl#1{_{\leq {#1}}}
\def\dg#1{_{> {#1}}}

\def\x{{\bf x}}
\def\y{{\bf y}}

\def\al{\alpha/\lambda}
\def\bl{\beta/\lambda}
\def\la{\lambda/\alpha}
\def\lb{\lambda/\beta}
\def\lm{\lambda/\mu}
\def\ln{\lambda/\nu}
\def\am{\alpha/\mu}
\def\bm{\beta/\mu}
\def\ma{\mu/\alpha} 
\def\mb{\mu/\beta}
\def\iff.{if and only if}
\def\rs{\stackrel{\rm\scriptscriptstyle RS}{\longleftrightarrow}}
\def\RS.{Robinson--Schensted correspondence}

\newdimen\psunit
\newdimen\psraisebox
\def\psraise(#1,#2)#3{\psraisebox=\psunit \multiply\psraisebox -#1
  \divide\psraisebox #2 \divide\psraisebox 2 \advance\psraisebox 3.2pt
\raisebox{\psraisebox}{#3}}
\psunit=12pt
\psset{unit=\psunit,linewidth=.4pt}
\newcount\ax \newcount\ay
\newcount\bx \newcount\by
\newcount\cx \newcount\cy
\newcount\dx \newcount\dy
\def\cell(#1,#2)[#3]{
\ax=#2 \ay=#1
\multiply\ay by-1
\bx=\ax \by=\ay
\cx=\ax \cy=\ay
\dx=\ax \dy=\ay
\advance\bx by-1
\advance\dy by1
\advance\cx by-1
\advance\cy by1
\pspolygon(\ax,\ay)(\bx,\by)(\cx,\cy)(\dx,\dy)(\ax,\ay)
\rput(\number\cx.5,
\ifnum\cy=0 -0.5\else\number\cy.5\fi){#3}
}

\section{Introduction}
In 2002, McKay, Morse and Wilf \cite{McKay2002} computed the large $n$ limit of
the probability that a random standard Young tableau of size $n$ contains a
fixed standard Young tableau. In this paper we obtain a $q$-analog of this
result. In order to state the McKay--Morse--Wilf theorem and our generalization,
we need to introduce several notations on permutations and standard Young
tableaux.

We denote by $\S_n$ the set of permutations of $[n]=\{1,2,\ldots,n\}$.  For a
permutation $\pi=\seq{\pi}{n}{}\in\S_n$, let $\pi\dl k$ (resp. $\pi\dg k$)
denote the permutation in $\S_k$ which is order-isomorphic to the subword of
$\pi$ consisting of the integers $i\leq k$ (resp. $i>k$). For example, if $\pi=5
1 3 6 9 7 4 2 8$, then $\pi\dl 4 = 1342$ and $\pi\dg4 = 12534$. Similarly, let
$\pi\ul k$ (resp. $\pi\ug k$) denote the permutation in $\S_k$ which is
order-isomorphic to $\seq{\pi}{k}{}$ (resp. $\pi_{k+1}\pi_{k+2}\cdots\pi_{n}$).
For example, if $\pi=5 1 3 6 9 7 4 2 8$, then $\pi\ul 4 = 3124$ and $\pi\ug4 =
53214$.  If $\sigma=\pi\dl k$ for some $k$, we say that $\pi$ {\em contains}
$\sigma$. Note that our permutation containment is different from the usual
pattern containment, see \cite{Bona}.

We denote by $\I_n$ the set of involutions in $\S_n$.  For a permutation
$\sigma\in \S_k$, let $\I_n(\sigma)$ denote the set of involutions
in $\I_n$ containing $\sigma$. In other words,
$$\I_n(\sigma)=\{\pi\in\I_n : \pi\dl k =\sigma\}.$$

A {\em partition} $\lambda =(\seq{\lambda}{r}{,})$ is a weakly decreasing
sequence of positive integers called \emph{parts}. The sum of the parts of a
partition $\lambda$ is denoted by $|\lambda|$. If $|\lambda|=n$, we say that
$\lambda$ is a partition of $n$, also written as $\lambda\vdash n$.  We denote
by $\emptyset$ the unique partition of $0$.  The {\em Ferrers diagram} of a
partition $\lambda=(\seq{\lambda}{r}{,})$ is the left justified array of squares
such that the $i$th row has $\lambda_i$ squares. For example, 
the Ferrers diagram of $(4,3,1)$ is the following:
\begin{center}
\begin{pspicture}(0,-3)(4,0)
    \cell(1,1)[] \cell(1,2)[] \cell(1,3)[] \cell(1,4)[]
    \cell(2,1)[] \cell(2,2)[] \cell(2,3)[] \cell(3,1)[]
 \end{pspicture}
\end{center}

We will identify a partition with its Ferrers diagram. Thus $\mu\subset\lambda$
means that the Ferrers diagram of $\mu$ is contained in the Ferrers diagram of
$\lambda$.  For two partitions $\lambda$ and $\mu$ with $\mu\subset\lambda$, the
{\em skew shape} $\lm$ is the set theoretic difference $\lambda-\mu$. If
$|\lambda|-|\mu|=n$, we write $\lm\vdash n$. For example, the skew shape
$\lm\vdash 4$ for $\lambda=(4,3,2)$ and $\mu=(3,2)$ is the following:
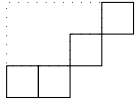
\begin{center}
\begin{pspicture}(0,-3) (4,0)
\cell(1,4)[] \cell(2,3)[] \cell(3,1)[] \cell(3,2)[]
\pspolygon[linestyle=dotted] (3,0) (3,-1)(2,-1) (2,-2) (0,-2) (0,0) 
\end{pspicture}
\end{center}
Occasionally we view a partition $\lambda$ as the skew shape
$\lambda/\emptyset$.

A {\em skew standard Young tableau} $T$, or a \emph{skew SYT} $T$ for short, is
a filling of a skew shape $\lm\vdash n$ with the integers in $[n]$ such that the
integers are increasing along rows and columns. In this case we say that the
\emph{shape} of $T$ is $\lm$, denoted by $\sh(T)=\lm$, and the \emph{size} of
$T$ is $n$. If the shape $T$ is a partition, $T$ is called a \emph{standard
  Young tableau} or \emph{SYT} for short. We denote by $\t_n$ the set of all SYTs of
size $n$.

For a SYT $T$, let $T_{\leq k}$ denote the SYT obtained by removing all squares
containing integers greater than $k$. Similarly let $T_{>k}$ denote the skew SYT
obtained by removing all squares containing integers less than or equal to $k$
and by decreasing all remaining integers by $k$.  For example, if
\begin{equation}
  \label{eq:5}
T= \psraise (3,1){\pspicture (0,-3) (4,0)
  \cell(1,1)[1] \cell(1,2)[2] \cell(1,3)[4] \cell(1,4)[7]
  \cell(2,1)[3] \cell(2,2)[5] \cell(2,3)[6] \cell(3,1)[8]
  \cell(3,2)[9] \endpspicture}  
\end{equation}
then $T\dl 5$ and $T\dg 5$ are the following:
\[
T\dl 5 = \psraise
(2,1){\pspicture (0,-2) (3,0) \cell(1,1)[1] \cell(1,2)[2]
  \cell(1,3)[4] \cell(2,1)[3] \cell(2,2)[5] \endpspicture} \qquad
T\dg 5= \psraise (3,1){\pspicture (0,-3) (4,0) \cell(1,4)[2] \cell(2,3)[1]
  \cell(3,1)[3] \cell(3,2)[4]\pspolygon[linestyle=dotted] (3,0) (3,-1)
  (2,-1) (2,-2) (0,-2) (0,0) \endpspicture}
\]

For two SYTs $U$ and $T$, if $U = T\dl k$ for some $k$, we say that $T$ {\em
  contains} $U$.  For a SYT $A$, we denote by $\t_n(A)$ the set of all SYTs of
size $n$ containing $A$.

We assume the reader is familiar with the \RS., see
\cite{Sagan2001,Stanley1999}.  If $\pi$ corresponds to the pair of SYTs $(P,Q)$
in the \RS., we write $\pi\rs(P,Q)$, and also $P(\pi)=P$ and $Q(\pi)=Q$. If
$\pi$ is an involution, we write $\pi\rs P$.

For given $A\in\t_m$, $\pi\in\I_n$ and $T\in\t_n$ with $\pi\rs T$, we can easily
see that $T$ contains $A$ \iff. $\pi$ contains $\sigma$ for some $\sigma\in\S_m$
with $P(\sigma)=A$. Thus the \RS. induces the following bijection:
\begin{equation}\label{eq:rs}
\mathrm{RS} : \t_n(A) \rightarrow \bigcup_{\sigma : P(\sigma)=A} \I_n(\sigma).  
\end{equation}

McKay, Morse and Wilf \cite{McKay2002} proved that for
$\sigma\in\S_m$, 
\begin{equation}\label{eq:mckay2}
\limit{n}\frac{|\I_{n}(\sigma)|}{|\I_n|} = 
\frac{1}{m!}.
\end{equation}
As a corollary, they obtained that the probability that a random SYT of size $n$
contains $A$ of shape $\alpha\vdash m$ tends to $f^{\alpha}/m!$ in the large $n$
limit, where $f^{\alpha}$ denotes the number of SYTs of shape $\alpha$.  In
other words,
\begin{equation}\label{eq:mckay1}
\limit{n} \frac{|\t_n(A)|}{|\t_n|}=
\frac{f^{\alpha}}{m!}.
\end{equation}

Jaggard~\cite{Jaggard2005} defined the $j$-set (see Section~\ref{s:perm})
$J(\pi)$ of a permutation $\pi$ and found the following exact formula for
$|\I_{n+m}(\sigma)|$ for $\sigma\in\S_m$:
\begin{equation}\label{eq:jaggard}
|\I_{n+m}(\sigma)|=\sum_{\substack{j\in J(\sigma)\\ k=n-m+j}} \binom{n}{k} t_k,
\end{equation}
where $t_k=|\I_k|$, the number of involutions in $\S_k$. Using
\eqref{eq:jaggard}, Jaggard~\cite{Jaggard2005} found another proof of
\eqref{eq:mckay2}.

In this paper we find $q$-analogs of \eqref{eq:mckay2}, \eqref{eq:mckay1} and
\eqref{eq:jaggard}. To state these, we need the following definitions.

For a permutation $\pi=\seq{\pi}{n}{}$, a {\em descent} of $\pi$ is an integer
$i\in[n-1]$ such that $\pi_{i}>\pi_{i+1}$.  For a skew SYT $T$ of size $n$, a
{\em descent} of $T$ is an integer $i\in[n-1]$ such that $i+1$ is in a row lower
than the row containing $i$.  Let $D(\pi)$ (resp. $D(T)$) denote the set of all
descents of $\pi$ (resp. $T$). Let $\maj(\pi)$ (resp. $\maj(T)$) denote the sum
of all descents of $\pi$ (resp. $T$). For example, if $\pi=5 1 3 6 9 7 4 2 8$,
we have $D(\pi)=\{1,5,6,7\}$ and $\maj(\pi)=19$. For the SYT $T$ in
\eqref{eq:5}, we have $D(T)=\{2,4,7\}$ and $\maj(T)=13$.

We define
$$\imaj(\pi)=\maj(\pi^{-1}),\qquad 
A_n(p,q)=\sum_{\pi\in\S_n} p^{\imaj(\pi)} q^{\maj(\pi)},$$
$$ t_n(q) = \sum_{\pi\in\I_n} q^{\maj(\pi)},\qquad f^{\lm}(q) =
\sum_{\sh(T)=\lm} q^{\maj(T)},$$
$$\qfactorial{n} =(1+q)(1+q+q^2)\cdots(1+q+\cdots+q^{n-1}),\qquad
\qbinom{n}{k} = \frac{\qfactorial{n}} {\qfactorial{k}
  \qfactorial{n-k}}.$$

Finally we can state our main results.
\begin{thm}\label{thm:main1}
  For $\sigma\in\S_m$, we have
$$\sum_{\pi\in \I_{n+m}(\sigma)} q^{\maj(\pi\ug m)} = 
\sum_{\substack{j\in J(\sigma)\\ k=n-m+j}} 
q^{\maj(\sigma\ug{j})} \qbinom{n}{k} t_k(q).$$
\end{thm}

For a real number $r>0$, we define $\rr=\min(r,r^{-1})$.

\begin{thm}\label{thm:main2}
  For $\sigma\in\S_m$ and a real number $q>0$, we have
$$\limit{n} \frac{\displaystyle\sum_{\pi\in\I_n(\sigma)}
  q^{\maj(\pi\ug m)} } {\displaystyle\sum_{\pi\in\I_n} q^{\maj(\pi\ug m)}}
=\frac{q^{\maj(\sigma)} + (1-\qq)C}{\qfactorial{m} + (1-\qq)D},$$
where $C$ and $D$ are polynomials in $q$ and $\qq$.
(We refer the reader to Theorem~\ref{thm:qlim1} for their exact form.)
\end{thm}

\begin{thm}\label{thm:main3}
  For a SYT $A$ of shape $\alpha\vdash m$ and a real number $q>0$, we have
$$\limit{n} \frac{\displaystyle\sum_{T\in\t_n(A)}q^{\maj(T\dg m)}}
{\displaystyle\sum_{T\in\t_n}q^{\maj(T\dg m)}} = \frac{f^{\alpha}(q) +
  (1-\qq)E}{\qfactorial{m}+ (1-\qq)D},$$ where $E$ and $D$ are
polynomials in $q$ and $\qq$.
(We refer the reader to Theorem~\ref{thm:m3} for their exact form.)
\end{thm}

We can obtain similar results by considering pairs of SYTs. For pairs $(P,Q)$
and $(A,B)$ of SYTs, we say that $(P,Q)$ {\em contains} $(A,B)$ if $P$ and $Q$
contain $A$ and $B$ respectively. We denote by $\t_n(A,B)$ the set of pairs
$(P,Q)$ of SYTs of the \emph{same shape} of size $n$ containing $(A,B)$.

Given $A\in\t_a$, $B\in\t_b$ and $\pi\rs (P,Q)$, it is easy to see that $(P,Q)$
contains $(A,B)$ \iff. $\pi\dl a = \sigma$ and $\pi\ul b=\tau$ for some
$\sigma\in\S_a$ and $\tau\in\S_b$ with $P(\sigma)=A$ and $Q(\tau)=B$. We denote
by $\sp$ the set of $\pi\in\S_n$ such that $\pi$ contains $\sigma$ and
$\pi^{-1}$ contains $\tau^{-1}$. In other words, for $\sigma\in\S_a$ and
$\tau\in\S_b$, 
$$\sp=\{\pi\in\S_n : \pi\dl a = \sigma, \pi\ul b=\tau \}.$$ Thus the
\RS. induces the following bijection:
\begin{equation}\label{eq:rs2}
RS : \t_n(A,B) \rightarrow 
\bigcup_{\substack{\sigma:P(\sigma)=A\\ \tau:Q(\tau)=B}} \sp.
\end{equation}

In Section~\ref{s:perm} we define the $j_2$-set $J(\sigma,\tau)$ of a pair
$(\sigma,\tau)$ of permutations, and prove the following $(p,q)$-analogs of
Theorems~\ref{thm:main1}, \ref{thm:main2} and \ref{thm:main3}.

\begin{thm}\label{thm:main1-1}
Let $a,b,n,m$ and $\ell$ be integers with $a+m=b+n=\ell$. Then, for 
$\sigma\in\S_a$ and $\tau\in\S_b$, we have
$$\sum_{\pi\in \spl} p^{\imaj(\pi \dg a)}q^{\maj(\pi\ug b)} =
\sum_{\substack{j\in J(\sigma,\tau)\\ k=n-a+j\\}}
p^{\imaj(\tau\dg{j})}q^{\maj(\sigma\ug{j})} \pbinom{m}{k}
\qbinom{n}{k} A_k(p,q).$$
In particular, when $p=q=1$ we have
$$|\spl|=\sum_{\substack{j\in J(\sigma,\tau)\\ k=n-a+j}} \binom{m}{k} \binom{n}{k} k!.$$
\end{thm}

\begin{thm}\label{thm:main2-1}
  For $\sigma\in\S_a$, $\tau\in\S_b$ and real numbers $p,q>0$, we have
$$\limit{n} \frac{\displaystyle\sum_{\pi\in\sp}p^{\imaj(\pi\dg a)}q^{\maj(\pi\ug b)}}
{\displaystyle\sum_{\pi\in\S_n}p^{\imaj(\pi\dg a)}q^{\maj(\pi\ug b)}} =
\frac{p^{\imaj(\tau)}q^{\maj(\sigma)} + (1-\pp)(1-\qq)C'}
{\pfactorial{b}\qfactorial{a} + (1-\pp)(1-\qq)D'},$$
where $C'$ and $D'$ are polynomials in $p,\pp,q$ and $\qq$.
(We refer the reader to Theorem~\ref{thm:m2-1} for their exact form.)
In particular, when $p=1$ or $p=q=1$ we have
$$\limit{n} \frac{\displaystyle\sum_{\pi\in\sp}q^{\maj(\pi\ug b)}}
     {\displaystyle\sum_{\pi\in\S_n}q^{\maj(\pi\ug b)}} =
     \frac{q^{\maj(\sigma)}}{b! \qfactorial{a}},$$
$$\limit{n}\frac{|\sp|}{|\S_n|} = \frac{1}{a!b!}.$$
\end{thm}

\begin{thm}\label{thm:main3-1}
  Let $A$, $B$ be SYTs of shape $\alpha\vdash a$, $\beta\vdash b$
  respectively. Then, for real numbers $p,q>0$, we have
$$\limit{n} \frac{\displaystyle\sum_{(P,Q)\in\t_n(A,B)}p^{\maj(P\dg a)}
  q^{\maj(Q\dg b)}}
{\displaystyle\sum_{(P,Q)\in\t_n(\emptyset,\emptyset)}p^{\maj(P\dg a)}q^{\maj(Q\dg
    b)}} = \frac{ f^{\beta}(p)f^{\alpha}(q) + (1-\pp)(1-\qq)E'}
{\pfactorial{b} \qfactorial{a} + (1-\pp)(1-\qq)D'},$$ where $E'$ and
$D'$ are polynomials in $p,\pp,q$ and $\qq$.  
(We refer the reader to Theorem~\ref{thm:m3-1} for their exact form.)
In particular, when $p=1$ or $p=q=1$ we have
 $$\limit{n} \frac{\displaystyle\sum_{(P,Q)\in\t_n(A,B)}q^{\maj(Q\dg b)}}
  {\displaystyle\sum_{(P,Q)\in\t_n(\emptyset,\emptyset)}q^{\maj(Q\dg b)}} =
    \frac{f^{\beta}f^{\alpha}(q)}{b!\qfactorial{a}},$$
$$\limit{n} \frac{|\t_n(A,B)|}{|\t_n(\emptyset,\emptyset)|} =
    \frac{f^{\alpha}f^{\beta}}{a!b!}.$$  
\end{thm}

The rest of this paper is organized as follows. In Section~\ref{s:perm}, we
define $j_2$-sets and prove Theorems~\ref{thm:main1}, \ref{thm:main2},
\ref{thm:main1-1} and \ref{thm:main2-1}.  In Section~\ref{s:tab}, we prove
Theorems~\ref{thm:main3} and \ref{thm:main3-1}.  In Section~\ref{sec:fcj}, we
find a simple method to determine whether a given set is a $j_2$-set. In
Section~\ref{sec:pf}, we prove a lemma which plays an important role in proving
the limit theorems.  In Section~\ref{sec:further-study}, we propose some open
problems.

\section{Permutation containment}\label{s:perm}
Recall the definitions of $\pi\ul k$, $\pi\ug k$, $\pi\dl k$ and $\pi\dg k$.
These are easy to remember using the following argument. We can consider a
permutation $\pi=\seq{\pi}{n}{}$ as a collection of {\em bi-letters}
$\biletter{i}{j}$ as follows:
$$\pi=\left\{ \biletter{1}{\pi_1}\quad \biletter{2}{\pi_2} \quad\cdots\quad 
  \biletter{n}{\pi_n}\right\}.$$ Then $\pi\dl k$ (resp. $\pi\dg k$, $\pi\ul k$
and $\pi\ug k$) is the permutation obtained from $\pi$ by taking the bi-letters
$\biletter{i}{j}$ with $j\leq k$ (resp.  $j>k$, $i\leq k$ and $i>k$) and by
relabeling them if necessary. It is easy to see that $(\pi\ul
k)^{-1}=(\pi^{-1})\dl k$.

In this section we prove Theorems~\ref{thm:main1}, \ref{thm:main2},
\ref{thm:main1-1} and \ref{thm:main2-1} which are related to permutation
containment. To do this, we start with a decomposition of permutation matrices.

For $\pi\in\S_n$, the {\em permutation matrix} $M(\pi)$ is the
$n\times n$ matrix whose $(i,j)$-entry is $1$ if $\pi_i=j$, and $0$
otherwise. For example, 
$$M(4132)=\mat{cccc}{0&0&0&1\\ 1&0&0&0\\ 0&0&1&0\\ 0&1&0&0}.$$

Consider a 0-1 matrix $M$ such that each row and column contains at most one 1.
There exists a unique permutation $\pi$ whose permutation matrix is obtained
from $M$ by removing the rows and columns consisting of zeroes. In this case, we
write $\pi\sim M$.  If $\pi\sim M$ and $\pi\sim N$, we also write $M\sim
N$. For example,
\begin{equation}
  \label{eq:3}
\mat{cccc}{ 0&1&0&0\\ 0&0&0&0\\ 1&0&0&0\\ 0&0&0&1} \sim
\mat{ccccc}{ 0&1&0&0&0\\ 0&0&0&0&0\\ 1&0&0&0&0\\ 0&0&0&1&0} \sim
213.
\end{equation}

For an $n\times m$ matrix $M$, we denote by $\row(M)$ (resp. $\col(M)$) the word
$\vr=r_1r_2\cdots r_n$ (resp.  $\vc=c_1c_2\cdots c_m$) of integers such that
$r_i$ (resp. $c_i$) is the sum of elements in the $i$th row (resp. column) of
$M$. For example, if $M$ is the second matrix in \eqref{eq:3}, then
$\row(M)=1011$ and $\col(M)=11010$.

Let $a,b,m,n$, and $\ell$ be fixed integers with $a+m=b+n=\ell$.  Consider
a permutation $\pi\in\S_\ell$.  We divide the permutation matrix of $\pi$ as
follows:
\[
M(\pi)=\bordermatrix{&a&m\cr b&M_{(1,1)}&M_{(1,2)}\cr
  n&M_{(2,1)}&M_{(2,2)}},  
\]
where the numbers outside the matrix
indicate the sizes of the block matrices. Note that
\[\pi\dl a \sim \mat{c}{M_{(1,1)}\\ M_{(2,1)}},\quad
\pi\dg a \sim \mat{c}{M_{(1,2)}\\ M_{(2,2)}},\]
\[\pi\ul b \sim \mat{cc}{M_{(1,1)} & M_{(1,2)}},\quad
\pi\ug b  \sim \mat{cc}{M_{(2,1)}& M_{(2,2)}}.\]

Assume that $M_{(1,1)}$ contains $j$ 1's. Then $M_{(1,2)}$,
$M_{(2,1)}$ and $M_{(2,2)}$ contain $b-j$, $a-j$ and $n-a+j$ 1's
respectively.
Then we define
$$\phi_{a,b}(\pi)=(\pi_{(1,1)},\pi_{(1,2)},\pi_{(2,1)},\pi_{(2,2)},
\vc_1,\vr_1,\vc_2,\vr_2),$$
where $\pi_{(r,s)}$ is the permutation
satisfying $\pi_{(r,s)} \sim M_{(r,s)}$ for $r=1,2$ and $s=1,2$, and
$\vc_1=\col(M_{(2,1)})$, $\vc_2=\col(M_{(2,2)})$,
$\vr_1=\col(M_{(1,2)})$ and $\vr_2=\col(M_{(2,2)})$.  
For example, if $\pi=7152436$ with
\[
M(\pi)=\left( \begin {tabular}{cc|ccccc}
  0&0&0&0&0&0&1\\ 1&0&0&0&0&0&0\\ 0&0&0&0&1&0&0\\ \hline 0&1&0&0&0
  &0&0\\ 0&0&0&1&0&0&0\\ 0&0&1&0&0&0&0 \\ 0&0&0&0&0&1&0\end {tabular}
  \right),
\]
then
\[
\phi_{2,3}(\pi)=(1,21,1,213,01,101,11010,0111).
\]

It is easy to see that the following is a bijection:
$$\phi_{a,b}:\S_\ell \rightarrow \bigcup_{\substack{0\leq j\leq
    a\\ k=n-a+j}} \S_{j}\times \S_{b-j}\times \S_{a-j}\times \S_k
\times \binom{[a]}{a-j} \times \binom{[b]}{b-j} \times \binom{[m]}{k}
\times \binom{[n]}{k},$$
where $\binom{[n]}{k}$ denotes the set of
words consisting of $k$ 1's and $n-k$ 0's.

For $\sigma\in\S_a$, $\tau\in\S_b$ and $\vr\in\binom{[a+b]}{b}$, the
{\em shuffle} $\sf(\sigma,\tau;\vr)$ is the permutation in $\S_{a+b}$
obtained from $\vr$ by replacing the $i$th 0 to $\sigma_i$ and the
$j$th 1 to $a+\tau_j$ for $1\leq i\leq a$ and $1\leq j\leq b$. For
example, $\sf(3142,231;0010110)=3164752$.  The following lemma is due
to Garsia and Gessel \cite{Garsia1979}.

\begin{lem}\label{thm:shuffle}
  For $\sigma\in\S_a$, $\tau\in\S_b$, we have
$$\sum_{\vr\in \binom{[a+b]}{b}} q^{\maj(\sf(\sigma,\tau;\vr))}=
q^{\maj(\sigma)+\maj(\tau)} \qbinom{a+b}{b}.$$
\end{lem}

Jaggard~\cite{Jaggard2005} defined the $j$-set as follows.  For a
permutation $\pi$, the {\em $j$-set} $J(\pi)$ of $\pi$ is defined to
be the set of integers $j\geq0$ such that $\pi\ul j$ is an
involution. In other words,
$$J(\pi) = \{ j: \pi\ul j = (\pi^{-1})\dl j \}.$$ 
Note that we always have $0\in J(\pi)$ because $\pi\ul 0 = (\pi^{-1})\dl 0$ is
the empty permutation. 

As we have seen in \eqref{eq:jaggard}, the $j$-set is useful to express the
number of elements in $\I_n(\sigma)$.  We define the $j_2$-set which is useful
to express the number of elements in $\sp$.

\begin{defn}
  The \emph{$j_2$-set} $J(\sigma,\tau)$ of a pair $(\sigma,\tau)$ of
  permutations is defined to be
$$J(\sigma, \tau) = \{ j : \sigma\ul j = \tau \dl j \}.$$
\end{defn}

Now we are ready to prove the following two theorems from which
Theorems~\ref{thm:main1} and \ref{thm:main1-1} follow. 

\begin{thm}\label{thm:permcont1}
For $\sigma\in \S_a$, we have
\begin{align*}
\sum_{\pi\in \I_{n+a}} q^{\maj(\pi\ug a)} &=
\sum_{\substack{0\leq j\leq a\\ k=n-a+j\\}}
t_j \binom{a}{j} \qfactorial{a-j} \qbinom{n}{k} t_k(q),\\
\sum_{\pi\in \I_{n+a}(\sigma)} q^{\maj(\pi\ug a)} &= 
\sum_{\substack{j\in J(\sigma)\\ k=n-a+j\\}}
q^{\maj(\sigma\ug{j})} \qbinom{n}{k} t_k(q).
\end{align*}
\end{thm}
\begin{proof}
We omit the proof since it is similar to the proof of the following theorem.
\end{proof}

\begin{thm}\label{thm:permcont2}
Let $a,b,n,m$ and $\ell$ be integers with $a+m=b+n=\ell$. Then, for 
$\sigma\in\S_a$ and $\tau\in\S_b$, we have
 \begin{align*}
\sum_{\pi\in \S_\ell} p^{\imaj(\pi \dg a)}q^{\maj(\pi\ug b)} &=
\sum_{\substack{0\leq j \leq a \\ k=n-a+j\\}}
j! \binom{a}{j} \binom{b}{j} \pfactorial{b-j} \qfactorial{a-j}
\pbinom{m}{k} \qbinom{n}{k} A_k(p,q),\\
\sum_{\pi\in \spl} p^{\imaj(\pi \dg a)}q^{\maj(\pi\ug b)} &=
\sum_{\substack{j\in J(\sigma,\tau)\\ k=n-a+j\\}}
p^{\imaj(\tau\dg{j})}q^{\maj(\sigma\ug{j})} \pbinom{m}{k}
\qbinom{n}{k} A_k(p,q).
\end{align*}
\end{thm}
\begin{proof}
  Fix an integer $j$ with $0\leq j\leq a$ and consider a permutation
  $\pi\in\S_\ell$ such that $\pi_{(1,1)}\in\S_j$, where
$$\phi_{a,b}(\pi)=(\pi_{(1,1)},\pi_{(1,2)},\pi_{(2,1)},\pi_{(2,2)},
\vc_1,\vr_1,\vc_2,\vr_2).$$ Then $(\pi\dg a)^{-1} =
\sf(\pi_{(1,2)}^{-1}, \pi_{(2,2)}^{-1}; \vc_2)$ and $\pi\ug b =
\sf(\pi_{(2,1)},\pi_{(2,2)}; \vr_2)$.  Thus
\begin{equation}\label{eq:1}
p^{\imaj(\pi \dg a)}q^{\maj(\pi\ug b)} = p^{\maj(\sf(\pi_{(1,2)}^{-1},
  \pi_{(2,2)}^{-1}; \vc_2))} q^{\maj(\sf(\pi_{(2,1)},\pi_{(2,2)}; \vr_2))}.
\end{equation}
Let $k=n-a+j$. By Lemma~\ref{thm:shuffle}, the sum of \eqref{eq:1} over
all $\pi_{(2,2)}\in \S_k$, $\vc_2\in\binom{[m]}k$ and
$\vr_2\in\binom{[n]}k$ equals
\begin{equation}\label{eq:2}
p^{\imaj(\pi_{(1,2)})}q^{\maj(\pi_{(2,1)})}\pbinom{m}{k}\qbinom{n}{k}A_k(p,q).
\end{equation}

Summing \eqref{eq:2} over all $j$, $\pi_{(1,1)}\in \S_j$, $\pi_{(1,2)}\in
\S_{b-j}$, $\pi_{(2,1)}\in \S_{a-j}$, $\vc_1\in\binom{[a]}{a-j}$ and
$\vr_1\in\binom{[b]}{b-j}$, and using the well known result $\sum_{\pi\in\S_n}
q^{\maj(\pi)}=\qfactorial{n}$, see \cite{Bona,Stanley1997}, we get the first
identity.

If $\pi\in\spl$, then $\sf(\pi_{(1,1)},\pi_{(1,2)};\vr_1)=\tau$ and
$\sf(\pi_{(1,1)}^{-1},\pi_{(2,1)}^{-1};\vc_1)=\sigma^{-1}$, which
implies $\pi_{(1,1)} = \sigma\ul j = \tau\dl j$, $\pi_{(1,2)}=\tau\dg
j$ and $\pi_{(2,1)}=\sigma\ug j$.  Thus we have $j\in J(\sigma,\tau)$,
and $j$ determines $\pi_{(1,1)}$, $\pi_{(1,2)}$, $\pi_{(2,1)}$,
$\vc_1$ and $\vr_1$. Then we get the second identity by summing
\eqref{eq:2} over all $j\in J(\sigma,\tau)$.
\end{proof}

To prove Theorems~\ref{thm:main2} and \ref{thm:main2-1} we need the following
lemma whose proof is given in Section~\ref{sec:pf}.  Recall that for a real number
$r>0$, we denote $\rr=\min(r,r^{-1})$.

\begin{lem}\label{thm:tlim}
For real numbers $p,q>0$, we have
$$\limit{n} \frac{\displaystyle \frac{t_{n+1}(q)}{\qfactorial{n+1}}}
{\displaystyle \frac{t_n(q)}{\qfactorial{n}}} = 1-\qq,\quad \limit{n}
\frac{\displaystyle \frac{A_{n+1}(p,q)}{\pfactorial{n+1}\qfactorial{n+1}}}
{\displaystyle \frac{A_n(p,q)}{\pfactorial{n}\qfactorial{n}}}
 = (1-\pp)(1-\qq).$$
\end{lem}

Now we prove the following theorem which implies Theorem~\ref{thm:main2}.

\begin{thm}\label{thm:qlim1}
  For $\sigma\in\S_m$ and a real number $q>0$, we have
$$\limit{n} \frac{\displaystyle\sum_{\pi\in\I_n(\sigma)}
  q^{\maj(\pi\ug m)} } {\displaystyle\sum_{\pi\in\I_n} q^{\maj(\pi\ug m)}}
=\frac{\displaystyle\sum_{j\in J(\sigma)} q^{\maj(\sigma\ug j)} \qbinom{m}{j}
\qfactorial{j}(1-\qq)^j}
{\displaystyle\sum_{j=0}^m \qfactorial{m} t_j \binom{m}{j} (1-\qq)^j}.$$
\end{thm}
\begin{proof}
By Theorem~\ref{thm:permcont1}, the left hand side is equal to
\begin{equation}\label{eq:pfeq1}
 \limit{n} \frac{\displaystyle \sum_{\pi\in\I_{n+m}(\sigma)} \frac{ q^{\maj(\pi\ug
      m)}}{\qfactorial n}} {\displaystyle \sum_{\pi\in\I_{n+m}} \frac{
    q^{\maj(\pi\ug m)}}{\qfactorial n}}= \limit{n} \frac{\displaystyle \sum_{j\in
    J(\sigma)} \frac{ q^{\maj(\sigma\ug j)}}{\qfactorial {m-j}}
  \frac{t_{n-m+j}(q)}{\qfactorial{n-m+j}}} {\displaystyle \sum_{j=0}^m t_j
  \binom{m}{j} \frac{t_{n-m+j}(q)}{\qfactorial{n-m+j}}}. 
\end{equation}
By Lemma~\ref{thm:tlim}, we have
$$\limit{n} \frac{\displaystyle \frac{t_{n-m+j}(q)}{\qfactorial{n-m+j}}}
{\displaystyle \frac{t_{n-m}(q)}{\qfactorial{n-m}}}= (1-\qq)^j.$$
Then, we get the theorem by dividing the numerator and the denominator of
the right hand side of \eqref{eq:pfeq1} by
$t_{n-m}(q)/\qfactorial{n-m}$, and by multiplying them by
$\qfactorial m$.
\end{proof}

Similarly, we can prove the following theorem, which implies
Theorem~\ref{thm:main2-1}.

\begin{thm}\label{thm:m2-1}
  For $\sigma\in\S_a$, $\tau\in\S_b$ and real numbers $p,q>0$, we have
  \begin{multline*}
\limit{n} \frac{\displaystyle\sum_{\pi\in\sp}p^{\imaj(\pi\dg a)}q^{\maj(\pi\ug
    b)}} {\displaystyle\sum_{\pi\in\S_n}p^{\imaj(\pi\dg a)}q^{\maj(\pi\ug b)}}
\\=\frac{\displaystyle\sum_{j\in J(\sigma,\tau)}  p^{\imaj(\tau\dg j)}q^{\maj(\sigma\ug j)}
\pbinom{b}{j} \qbinom{a}{j} \pfactorial j \qfactorial j
(1-\pp)^j(1-\qq)^j}{\displaystyle\sum_{j=0}^a \pfactorial b \qfactorial a j! 
\binom{a}{j} \binom{b}{j}(1-\pp)^j(1-\qq)^j}.
 \end{multline*}
\end{thm}

Note that when $p=q=1$ in Theorem~\ref{thm:m2-1} we have
\begin{equation}
  \label{eq:4}
\limit{\ell}\frac{|\sp|}{|\S_\ell|} = \frac{1}{a!b!},
\end{equation}
which means that the probability that a random permutation $\pi\in\S_\ell$
contains $\sigma$ and $\pi^{-1}$ contains $\tau^{-1}$ tends to $1/a!b!$
as $\ell$ approaches infinity.  This statement fits well with the following
intuition. If we divide $M(\pi)$ as
\[
M(\pi)=\bordermatrix{&a&m\cr b&M_{(1,1)}&M_{(1,2)}\cr n&M_{(2,1)}&M_{(2,2)}},  
\]
where $m=\ell-a$ and $n=\ell-b$, then the probability that $M_{(1,1)}$ contains
a $1$ will be very low when $\ell$ is very large. Thus the probability that 
$M_{(2,1)}\sim\sigma$ and $M_{(1,2)}\sim\tau$ will be very close to $1/a!$ and
$1/b!$ respectively, which is consistent with \eqref{eq:4}.

\section{Tableau containment}\label{s:tab}

Jaggard \cite{Jaggard2005} proved that for a SYT $A$ of shape $\alpha$
and an integer $j$,
\begin{equation}\label{eq:jag}
\#\{\sigma:P(\sigma)=A, j\in J(\sigma)\} = \sum_{\mu\vdash j} f^{\am}.
\end{equation}

We have the following analog of \eqref{eq:jag}: given SYTs $A$ and $B$
of shape $\alpha$ and $\beta$ respectively and an integer $j$, we have
\begin{equation}
  \label{eq:jag-j2}
\#\{(\sigma,\tau): P(\sigma)=A, Q(\tau)=B, j\in J(\sigma,\tau)\}
 = \sum_{\mu\vdash j} f^{\bm} f^{\am}.
\end{equation}

Using \eqref{eq:jag}, Jaggard \cite{Jaggard2005} gave another proof of the following
formula due to Sagan and Stanley \cite{Sagan1990}:
\begin{equation}\label{eq:sag1}
\sum_{\la\vdash n} f^{\la} = \sum_{k\geq0} \binom{n}{k} t_k
\sum_{\am\vdash n-k} f^{\am}.  
\end{equation}

Similarly using \eqref{eq:jag-j2} we can prove the following formula also due to Sagan and Stanley
\cite{Sagan1990}:
\begin{equation}\label{eq:sag2}
\sum_{\substack{\la \vdash m\\ \lb\vdash n}} f^{\la}f^{\lb} =
\sum_{k\geq0} \binom{m}{k} \binom{n}{k} k! \sum_{\substack{\bm \vdash
    m-k\\ \am\vdash n-k}} f^{\bm} f^{\am}.  
\end{equation}

In this section we generalize these results and prove Theorems~\ref{thm:main3}
and \ref{thm:main3-1}.

We will use the following well known fact: for given $\pi\rs (P,Q)$ we have
$D(\pi^{-1})=D(P)$ and $D(\pi)=D(Q)$, see \cite{Stanley1999}.  Thus
$q^{\maj((\pi^{-1})\ug a)}= q^{\maj(P\dg a)}$ and $q^{\maj(\pi\ug b)}=
q^{\maj(Q\dg b)}$ for all nonnegative integers $a$ and $b$.

We start with the following two lemmas which are
respectively a $q$-analog of \eqref{eq:jag} and a $(p,q)$-analog of
\eqref{eq:jag-j2}.

\begin{lem}\label{thm:permtotab}
  For a SYT $A$ of shape $\alpha$ and an integer $j$, we have
\[
\sum_{\sigma:\left\{\substack{P(\sigma)=A \\j\in J(\sigma)}\right.} 
q^{\maj(\sigma\ug j)} = \sum_{\mu\vdash j} f^{\am}(q).
\]
\end{lem}
\begin{proof}
Since it is similar to the proof of the following lemma, we omit it.
\end{proof}

\begin{lem}
For SYTs $A$ and $B$ of shape $\alpha$ and $\beta$ respectively and
an integer $j$, we have
$$\sum_{(\sigma,\tau): \left\{\substack{P(\sigma)=A\\ Q(\tau)=B \\j\in J(\sigma,\tau)}\right.}
p^{\imaj(\tau\dg j)} q^{\maj(\sigma\ug j)} = \sum_{\mu\vdash j} f^{\bm}(p) f^{\am}(q).$$
\end{lem}
\begin{proof}
Let $X=\{(\sigma,\tau):P(\sigma)=A,Q(\tau)=B,j\in J(\sigma,\tau)\}$
and $Y=\{(U,V): \mu\vdash j, \sh(U)=\bm, \sh(V)=\am\}$. It is sufficient
to find a bijection $\xi:X \rightarrow Y$ such that if
$\xi(\sigma,\tau)=(U,V)$, then $\imaj(\tau\dg j)=\maj(U)$ and
$\maj(\sigma\ug j) =\maj(V)$.

We define $\xi(\sigma,\tau)=(U,V)$ by $U=P(\tau)\dg j$ and
$V=Q(\sigma)\dg j$. Then we have $\imaj(\tau\dg j)=\maj(U)$ and
$\maj(\sigma\ug j) =\maj(V)$.  

To prove $\xi$ is a bijection, it is sufficient to show that for
$(U,V)\in Y$ there exists a unique pair $(\sigma,\tau)\in X$
satisfying $\xi(\sigma,\tau)=(U,V)$.  Let $\alpha\vdash a$. Since
$P(\sigma)=A$, $Q(\sigma)\dg j=V$, by reversing the insertion
algorithm $a-j$ times, we can find $\sigma_a,
\sigma_{a-1},\ldots,\sigma_{j+1}$ and $P(\sigma\ul j)$. Since
$P(\tau)\dl j = P(\tau\dl j)=P(\sigma\ul j)$ and $P(\tau)\dg j=U$, we
can determine $P(\tau)$. Thus we get $\tau$, from which we can
determine $\sigma\ul j = \tau\dl j$. Thus we get $\sigma$, and there
is a unique pair $(\sigma,\tau)\in X$ with $\xi(\sigma,\tau)=(U,V)$.
\end{proof}

Now we have the following $q$-analog of \eqref{eq:sag1} and $(p,q)$-analog of
\eqref{eq:sag2}.

\begin{thm}\label{thm:majgen}
  For a partition $\alpha$ and an integer $n$, we have
$$\sum_{\la\vdash n} f^{\la}(q) = \sum_{k\geq0} \qbinom{n}{k} t_k(q)
\sum_{\am\vdash n-k} f^{\am}(q).$$
\end{thm}
\begin{proof}
  Let $\alpha\vdash a$ and let $A$ be a SYT of shape $\alpha$. Then the left
  hand side is equal to
  \begin{align*}
    \sum_{T\in\t_{n+a}(A)} q^{\maj(T\dg a)} &= \sum_{\sigma:P(\sigma)=A}
    \sum_{\pi\in \I_{n+a}(\sigma)} q^{\maj(\pi\ug a)} \\
    &= \sum_{\sigma:P(\sigma)=A} \sum_{\substack{j\in J(\sigma)\\ k=n-a+j}}
       q^{\maj(\sigma\ug j)} \qbinom{n}{k} t_k(q) 
       & \mbox{(by Theorem~\ref{thm:permcont1})} \\
    &= \sum_{\substack{0\leq j\leq a\\ k=n-a+j}} \qbinom{n}{k} t_k(q)
       \sum_{\sigma : \left\{ \substack{P(\sigma)=A\\ j\in J(\sigma)} \right.} 
       q^{\maj(\sigma\ug j)}\\
    &= \sum_{\substack{0\leq j\leq a\\ k=n-a+j}} \qbinom{n}{k} t_k(q)
       \sum_{\mu\vdash j} f^{\am}(q)
       &\mbox{(by Lemma~\ref{thm:permtotab})}\\
    &= \sum_{k\geq0} \qbinom{n}{k} t_k(q) \sum_{\am\vdash n-k} f^{\am}(q).
  \end{align*}
\end{proof}

\begin{thm}\label{thm:majgen1}
  For partitions $\alpha$, $\beta$ and integers $m,n$, we have
$$\sum_{\substack{\la \vdash m\\ \lb\vdash n}} f^{\la}(p) f^{\lb}(q) =
\sum_{k\geq0} \pbinom{m}{k} \qbinom{n}{k} A_k(p,q) \sum_{\substack{\bm
    \vdash m-k\\ \am\vdash n-k}} f^{\bm}(p) f^{\am}(q).$$
\end{thm}
\begin{proof}
  This can be done similarly as in the proof of Theorem~\ref{thm:majgen}.
\end{proof}

We note that Theorems~\ref{thm:majgen} and \ref{thm:majgen1} can also
be proved using the following identities of skew Schur functions:
$$\sum_{\lambda} s_{\la}(\x)= \sum_{\mu} s_{\am}(\x)
\prod_{i}(1-x_i)^{-1} \prod_{i<j}(1-x_ix_j)^{-1},$$
$$\sum_{\lambda} s_{\la}(\x) s_{\lb}(\y) = \prod_{i,j}(1-x_i y_j)^{-1}
\sum_{\mu} s_{\bm}(\x) s_{\am}(\y).$$

Now we can prove the following theorem, which implies
Theorem~\ref{thm:main3}.
\begin{thm}\label{thm:m3}
  For a SYT $A$ of shape $\alpha\vdash m$ and a real number $q>0$, we have
$$\limit{n} \frac{\displaystyle\sum_{T\in\t_n(A)} q^{\maj(T\dg m)} }
{\displaystyle\sum_{T\in\t_n} q^{\maj(T\dg m)}} = 
\frac{\displaystyle \sum_{j=0}^m
  \qbinom{m}{j} \qfactorial{j}(1-\qq)^j \sum_{\mu\vdash j} f^{\am}(q)}
{\displaystyle\sum_{j=0}^m \qfactorial{m} t_j \binom{m}{j} (1-\qq)^j}.$$
\end{thm}
\begin{proof}
By the \RS. \eqref{eq:rs}, we have
$$\limit n \frac{\displaystyle\sum_{T\in\t_n(A)} q^{\maj(T\dg m)}}
{\displaystyle\sum_{T\in\t_n} q^{\maj(T\dg m)}} =\limit n
\frac{\displaystyle \sum_{\sigma:P(\sigma)=A} \sum_{\pi\in \I_{n}(\sigma)} q^{\maj(\pi\ug m)}}
{\displaystyle\sum_{\pi\in\I_n}q^{\maj(\pi\ug m)}}.$$
By Theorem~\ref{thm:qlim1}, the above limit is equal to
\[\frac{\displaystyle\sum_{\sigma:P(\sigma)=A}
  \sum_{j\in J(\sigma)} q^{\maj(\sigma\ug j)}
  \qbinom{m}{j}\qfactorial{j}(1-\qq)^j} {\displaystyle \sum_{j=0}^m \qfactorial{m}
  t_j \binom{m}{j} (1-\qq)^j}
=\frac{\displaystyle\sum_{j=0}^m \qbinom{m}{j}\qfactorial{j}(1-\qq)^j
\sum_{\sigma : \left\{ \substack{P(\sigma)=A\\ j\in J(\sigma)} \right.} 
q^{\maj(\sigma\ug j)}}{\displaystyle \sum_{j=0}^m \qfactorial{m}
  t_j \binom{m}{j} (1-\qq)^j}.\]
By Lemma~\ref{thm:permtotab}, we are done.
\end{proof}

Similarly, we get the following, which implies Theorem~\ref{thm:main3-1}.

\begin{thm}\label{thm:m3-1}
Let $A$ and $B$ be SYTs of shape $\alpha\vdash a$ and $\beta\vdash b$
respectively. Then for real number $p,q>0$, we have
\begin{multline*}
\limit{n} \frac{\displaystyle\sum_{(P,Q)\in\t_n(A,B)}p^{\maj(P\dg a)}q^{\maj(Q\dg
    b)}} {\displaystyle\sum_{(P,Q)\in\t_n(\emptyset,\emptyset)}p^{\maj(P\dg
    a)}q^{\maj(Q\dg b)}} \\
= \frac{\displaystyle\sum_{j=0}^a 
\pbinom{b}{j} \qbinom{a}{j} \pfactorial j \qfactorial j
(1-\pp)^j(1-\qq)^j
\sum_{\mu\vdash j}  f^{\bm}(p) f^{\am}(q)}
{\displaystyle\sum_{j=0}^a \pfactorial b \qfactorial a j! 
\binom{a}{j} \binom{b}{j}(1-\pp)^j(1-\qq)^j}.
\end{multline*}
\end{thm}

Let us consider the special cases $q=1$ and $p=q=1$ of Theorems~\ref{thm:m3} and
\ref{thm:m3-1}, which are
\begin{equation}
  \label{eq:6}
\limit{n} \frac{|\t_n(A)|}{|\t_n|}=
\frac{f^{\alpha}}{a!}.
\end{equation}
\begin{equation}
  \label{eq:7}
\limit{n} \frac{|\t_n(A,B)|}{|\t_n(\emptyset,\emptyset)|} =
    \frac{f^{\alpha}f^{\beta}}{a!b!}.  
\end{equation}

As mentioned in the introduction, \eqref{eq:6} is first proved by McKay et
al. \cite{McKay2002}. To our knowledge \eqref{eq:7} is new.  For the rest of
this section we give another proofs of \eqref{eq:6} and \eqref{eq:7} using
\eqref{eq:sag1} and \eqref{eq:sag2}.

Using the well known asymptotic behavior of $t_n \sim
\frac{1}{\sqrt2}n^{n/2} \exp\left(-\frac{n}{2}+\sqrt n
-\frac{1}{4}\right)$, we can easily see that 
\[ \limit{n} \frac{t_n}{t_{n+1}} = 0.\]
Thus
\[ \limit{n} \frac{n t_{n-1}}{t_{n+1}} = \limit{n} \frac{t_{n+1}-t_n}{t_{n+1}} =1.\]
Using induction one can easily prove that
\begin{equation}
  \label{eq:8}
 \limit{n} \frac{n^a t_{n-a}}{t_{n+a}} = \limit{n} \frac{n^a t_{n}}{t_{n+2a}} = 1.
\end{equation}

By \eqref{eq:sag1}, we have
\[
\limit{n}\frac{|\I_{n}(\sigma)|}{|\I_n|} = 
\limit{n}\frac{\displaystyle\sum_{\la\vdash n} f^{\la}}{t_{n+a}} = 
\limit{n}\frac{\displaystyle\sum_{k\geq0} \binom{n}{k} t_k
\sum_{\am\vdash n-k} f^{\am}}{t_{n+a}}.
\]
Since $\am\vdash n-k$, we only need to consider $k$ with $n-a\leq k\leq n$.
Using \eqref{eq:8}, it is not difficult to see that for $k>n-a$, we have
\[ \limit{n}\frac{\binom{n}{k} t_k}{t_{n+a}} = 0, \]
and
\[ \limit{n}\frac{\binom{n}{n-a} t_{n-a}}{t_{n+a}} = \frac{1}{a!}. \]
This gives another proof of \eqref{eq:6}.

The same argument can be applied to \eqref{eq:7}. More precisely, by
\eqref{eq:sag2},
\[
\limit{\ell} \frac{|\t_\ell(A,B)|}{|\t_\ell(\emptyset,\emptyset)|} =
\limit{\ell} \frac{\displaystyle\sum_{\substack{\la \vdash m\\ \lb\vdash n}} f^{\la}f^{\lb}}
{\ell!} =\limit{\ell} 
\frac{\displaystyle
\sum_{k\geq0} \binom{m}{k} \binom{n}{k} k! \sum_{\substack{\bm \vdash
    m-k\\ \am\vdash n-k}} f^{\bm} f^{\am}} {\ell!},
\]
where $\ell=a+m = b+n$, $\sh(A)=\alpha\vdash a$ and $\sh(B)=\beta\vdash b$.  By
the condition $\bm\vdash m-k$, we only need to consider the integers $k$ with $k\geq
m-b = n-a=\ell-a-b$. It is easy to see that for $k>m-b$, we have
\[ \limit{\ell} \frac{\binom{m}{k} \binom{n}{k} k!} {\ell!} = 0, \]
and
\[ \limit{\ell}  \frac{\binom{m}{m-b} \binom{n}{n-a} (\ell-a-b)!} {\ell!} =
\frac1{a!b!}. \]
This gives another proof of \eqref{eq:7}.

We note that in fact the above two proofs are special cases of the proofs of
Theorems~\ref{thm:m3} and \ref{thm:m3-1}.

\section{Criterion for $j_2$-sets}\label{sec:fcj}

In this section we find a criterion for a set to be a $j_2$-set.  First, we
review some previous results on $j$-sets. Throughout this section we
assume that $n$, $m$ and $k$ are positive integers.

Kim and Kim \cite{Kim2007} proved the following theorem.

\begin{thm}
\label{thm:kimkim}
Let $J$ be a $j$-set with largest element $m\geq 2$. Then, for
$n>m$, $J\cup \{n\}$ is a $j$-set if and only if $n=m+1$ or $n-m\geq
m-\max(J\cap [m-2])$.
\end{thm}

Using Theorem~\ref{thm:kimkim} we can easily determine whether a given set is
a $j$-set. To do this we need the following definition introduced by Corteel and
Lovejoy \cite{Corteel2004}.  An \emph{overpartition} is a weakly decreasing
sequence of positive integers in which the last occurrence of an integer may be
overlined.

For a set $S=\{s_0,s_1,s_2,\ldots,s_n\}$ of nonnegative integers with
$s_0<s_1<\cdots<s_n$,  we define $\Delta(S)$ to be the sequence
$(a_1,a_2,\ldots,a_n)$ where $a_i = s_{n-i+1}-s_{n-i}$ for $i\in[n]$. 
For example, if 
\begin{equation}
  \label{eq:9}
S=\{0,1,2,3,5,6,9,13,17,18,19,20,22\},
\end{equation}
then
\begin{equation}
  \label{eq:10}
\Delta(S)= (2,1,1,1,4,4,3,1,2,1,1,1).  
\end{equation}
And then we define $\ov\Delta(S)$ to be the sequence obtained from
$\Delta(S)=(a_1,a_2,\ldots,a_n)$ by applying the following algorithm.

First we set $k=1$. Find the smallest index $i$ such that $a_i=1$ and $k\leq
i<n$. If we have no such $i$, then we finish the algorithm. Otherwise we replace
$a_i$ and $a_{i+1}$ with $\ov{a_{i+1}+1}$. Set $k=i+2$ and repeat this
process.

For example, for $\Delta(S)$ in \eqref{eq:10}, we have
\begin{equation}
  \label{eq:11}
\ov{\Delta}(S)= (2,\ov2,\ov5,4,3,\ov3,\ov2,1).  
\end{equation}

Assume that we have $\ov\Delta(S)=(b_1,\ldots,b_m)$. Let ${\bf i}=\{i_1, i_2,
\ldots ,i_k\}$ be the set of integers with $1\leq i_1<i_2<\cdots <i_k\leq m$
such that $b_j=\ov2$ if and only if $j\in {\bf i}$. We define $\psi(S)$ to be
the sequence $(\ss_1,\ss_2,\ldots,\ss_{k+1})$, where $\ss_j = (b_{i_{j-1}+1},
b_{i_{j-1}+2}, \ldots, b_{i_j})$ for $j\in[k+1]$, and $i_0=0$ and $i_{k+1}=n$.
 
For example, for $S$ in \eqref{eq:9}, we have
\begin{equation}
  \label{eq:12}
  \psi(S) = ((2,\ov2),(\ov5,4,3,\ov3,\ov2),(1)).  
\end{equation}

With the above observation Kim and Kim \cite{Kim2007} proved the following.

\begin{cor}\label{thm:1}
  Let $S$ be a set of nonnegative integers with
  $\psi(S)=(\ss_1,\ldots,\ss_{k+1})$. Then $S$ is a $j$-set if and only if $0\in
  S$, $\ss_i$ is an overpartition ending with $\ov2$ for all $i\in[k]$ and
  $\ss_{k+1}=(1)$ or $\ss_{k+1}=\emptyset$.
\end{cor}

Hence, the set $S$ in \eqref{eq:9} is a $j$-set because $\psi(S)$ in
\eqref{eq:12} satisfies the condition in Corollary~\ref{thm:1}.

Note that since $J(\pi)=J(\pi,\pi^{-1})$, every $j$-set is also a
$j_2$-set. However, the converse is not true. For example, the $j_2$-set
$J(312,312)=\{0,1,3\}$ is not a $j$-set.

We have the following theorem analogous to Theorem~\ref{thm:kimkim}.

\begin{thm}\label{thm:j2}
Let $J$ be a $j_2$-set such that the two largest elements of $J$ are
$n-k$ and $n$. Then $J\cup\{n+m\}$ is a $j_2$-set \iff. $m=1$ or
$m\geq k$. 
\end{thm}

Before proving Theorem~\ref{thm:j2} we will see how to determine whether a
given set is a $j_2$-set.

Consider a set $S$ of integers with $\Delta(S)=(a_1,a_2,\ldots,a_n)$.  Let ${\bf
  i}=\{i_1, i_2, \ldots ,i_k\}$ be the set of integers with $1\leq
i_1<i_2<\cdots <i_k\leq n$ such that $a_j=1$ if and only if $j\in {\bf i}$. We
define $\psi_2(S)$ to be the sequence $(\ss_1,\ss_2,\ldots,\ss_k)$, where $\ss_j
= (a_{i_{j-1}+1}, a_{i_{j-1}+2}, \ldots, a_{i_j})$ for $j\in[k]$ and $i_0=0$.
For example, if
\begin{equation}
  \label{eq:13}
S=\{0,1,3,6,7,8,12,13,14,15,17\},  
\end{equation}
then
\begin{align*}
  \Delta(S) &=(2, 1,1,1,4,1,1, 3,2,1),\\
\psi_2(S) &= ((2,1),(1),(1),(4,1),(1),(3,2,1)).  
\end{align*}

It is easy to see that Theorem~\ref{thm:j2} implies the following.

\begin{cor}\label{thm:cri}
  Let $S$ be a set of nonnegative integers with
  $\psi_2(S)=(\ss_1,\ss_2,\ldots,\ss_k)$. Then $S$ is a $j_2$-set if and only if
  $0\in S$ and $\ss_i$ is a partition with exactly one part equal to $1$ for all
  $i\in [k]$.
\end{cor}

By Corollary~\ref{thm:cri}, the set $S$ in \eqref{eq:13} is a $j_2$-set. Using
Corollary~\ref{thm:cri}, we can easily get the following generating function for the
number $j_2(n)$ of $j_2$-sets with largest element $n$.

\begin{cor}
We have
\begin{align*}
\sum_{n\geq 0} j_2(n)x^n &
=\frac{1}{\displaystyle 1-x\prod_{i\geq2}\frac{1}{1-x^i}}\\
&=1 + x + x^2 + 2x^3 + 4x^4 + 8x^5 + 15 x^6 + 29x^7 + 55x^8 +105x^9\\
&\quad +200x^{10}+381x^{11}+725x^{12}+1381x^{13}+2629x^{14}+5005x^{15}+\cdots.
\end{align*}
\end{cor}

For the rest of this section we prove Theorem~\ref{thm:j2}. Our proof is similar to,
but simpler than, the proof of Theorem~\ref{thm:kimkim} in \cite{Kim2007}.  

Note that if $J$ is a $j_2$-set, then $J\cap [k]$ is also a $j_2$-set for all
integers $k$.

\begin{prop}
Let $J$ be a $j_2$-set with largest element $n$.
Then there is a permutation in $\S_n$ such that $J(\pi,\pi)=J$.
\end{prop}
\begin{proof}
  Let $(\sigma,\tau)$ be a pair with $J=J(\sigma,\tau)$.  Since $\sigma\ul
  n=\tau\dl n$, if we set $\pi=\sigma\ul n$, we have $J(\pi,\pi)=J$.
\end{proof}

Recall that, for 0-1 matrices $M$ and $N$, we write $M\sim N$ if the
matrices obtained from $M$ and $N$ by removing rows and columns
consisting of zeroes are the same.

To prove Theorem~\ref{thm:j2} we need the following four lemmas.

\begin{lem}\label{thm:2}
  Let $J$ be a $j_2$-set such that the three largest elements of $J$ are $n-k$,
  $n$ and $n+m$ with $m\geq2$. Then $m\geq k$.
\end{lem}
\begin{proof}
Consider  $\pi\in\S_{n+m}$ with $J(\pi,\pi)=J$. We divide
  $M(\pi)$ as follows:
  \[
M(\pi)=\bordermatrix{&n&m\cr n&A&B\cr m&C&D}.
\]
Since $n\in J(\pi,\pi)$, we have $\mat{cc}{A&B} \sim \mat{c}{A\\C}$.  Let
$\sigma=\pi\ul n=\pi\dl n$.  Assume that $B$ has $s$ nonzero entries. Then $C$
also has $s$ nonzero entries.  If $s=0$, then we get $n+1\in J(\pi,\pi)$ because
$\pi\ul{n+1}\sim \mat{cc}{A&0\\0&1} \sim \pi\dl{n+1}$. But this is a
contradiction to $m\geq2$.  Thus $s\geq 1$. Since $\sigma\ul{n-s}\sim A \sim
\sigma\dl{n-s}$, we get $\sigma\ul{n-s}=\sigma\dl{n-s}$. Thus
$$\pi\ul{n-s}=(\pi\ul n)\ul{n-s} = \sigma\ul{n-s} =
\sigma\dl{n-s}=(\pi\dl n)\dl{n-s}=\pi\dl{n-s}$$ and we get $n-s\in
J(\pi,\pi)$. Since $n-k$ is the largest element in $J\cap [n-1]$, we
get $n-s\leq n-k$.  Since $B$ has at most $m$ nonzero entries, we get
$k\leq s\leq m$.
\end{proof}

\begin{lem}\label{thm:3}
Let $J$ be a $j_2$-set such that the two largest elements of $J$ are
$n-k$ and $n$. Then $J\cup\{n+k\}$ is a $j_2$-set.
\end{lem}
\begin{proof}
  If $k=1$, then it is clear. Assume $k\geq 2$.  Consider $\sigma\in\S_n$ with
  $J(\sigma,\sigma)=J$. Let $\pi$ be the unique permutation in $\S_{n+k}$
  satisfying
\begin{equation}\label{eq:pi}
M(\pi)=\bordermatrix{&n&k\cr n&A&C\cr k&B&\bf 0},
\end{equation}
where $A,B$ and $C$ are the matrices of size $n\times n$, $k\times n$ and
$n\times k$ respectively such that
$M(\sigma)\sim\mat{c}{A\\B}\sim\mat{cc}{A&C}$. The condition $n-k\in
J(\sigma,\sigma)$ guarantees that there is a unique $\pi$. Since $J(\pi,\pi)\cap
[n]=J$ and $n+k\in J(\pi,\pi)$, it is sufficient to show that $n+s\notin
J(\pi,\pi)$ for all $1\leq s<k$. Suppose $n+s\in J(\pi,\pi)$ for some $1\leq
s<k$. Then we have
\begin{equation}\label{eq:mat}
\mat{cc}{A&C\\B'& \bf 0}\sim \mat{cc}{A&C'\\B& \bf 0},
\end{equation}
where $B'$ (resp. $C'$) is the matrix consisting of the first $s$ rows
of $B$ (resp. columns of $C$). Removing the last $k-s$ nonzero rows
and columns of the matrices in both sides of \eqref{eq:mat}, we get
$\sigma\ul{n-k+s}=\sigma\dl{n-k+s}$, i.e., $n-k+s\in
J(\sigma,\sigma)=J$, which is a contradiction to the assumption that
$n-k$ and $n$ are the two largest elements of $J$.
\end{proof}

\begin{lem}\label{thm:4}
Let $J$ be a $j_2$-set such that the two largest elements of $J$ are
$n-k$ and $n$. If $k\geq 2$, then $(J\setminus\{n\})\cup\{n+1\}$ is a
$j_2$-set.
\end{lem}
\begin{proof}
Consider $\sigma=\seq{\sigma}{n}{}\in\S_n$ with
$J(\sigma,\sigma)=J$. Since $n-1\notin J(\sigma,\sigma)$, we have
$\sigma_n\ne n$.  Let $\pi\in \S_{n+1}$ be the permutation such that
$$\pi_i = \left\{
\begin{array}{ll}
  \sigma_i & \mbox{if $i<n$ and $\sigma_i<n$,}\\
  n+1 & \mbox{if $i<n$ and $\sigma_i=n$,}\\
  n & \mbox{if $i=n$,}\\
  \sigma_n & \mbox{if $i=n+1$.}
\end{array}\right.$$
Then $M(\sigma)$ and $M(\pi)$ are decomposed as follows:
$$M(\sigma)=\begin{tabular}{|ccc|c|} \hline
  &  &  &  \\ 
  & $A$ &  & $C$ \\
  &  &  &  \\  \hline
  & $B$ &  & $0$ \\ \hline
\end{tabular}\:, \quad
M(\pi)=\begin{tabular}{|ccc|c|c|} \hline
  &  &  & $0$ &  \\ 
  & $A$ &  & $\vdots$ & $C$ \\
  &  &  & $0$ &  \\  \hline
 0 & $\cdots$ & 0 & 1 & 0 \\ \hline
  & $B$ &  & 0 & 0 \\ \hline 
\end{tabular}\:,$$
where $A$, $B$ and $C$ are $(n-1)\times(n-1)$, $1\times (n-1)$ and
$(n-1)\times 1$ matrices respectively.  It is not difficult to see
that $J(\pi,\pi)=(J\setminus\{n\})\cup\{n+1\}$.
\end{proof}

\begin{lem}\label{thm:5}
Let $J$ be a $j_2$-set such that the two largest elements of $J$ are
$n-1$ and $n$. Then $J\cup\{n+k\}$ is a $j_2$-set for all positive
integers $k$.
\end{lem}
\begin{proof}
  It is clear if $k=1$. Assume $k\geq 2$.  Consider $\sigma\in \S_{n-1}$ with
  $J(\sigma,\sigma)=J\setminus\{n\}$.  Let $\pi$ be the permutation in
  $\S_{n+k}$ such that $M(\pi)=\mat{cc}{M(\sigma)&{\bf 0}\\{\bf 0}&A}$, where
  $A=\mat{cc}{{\bf 0} & I_{k-1} \\ 1 & {\bf 0}}$ and $I_{k-1}$ is the
  $(k-1)\times (k-1)$ identity matrix. Then $J(\pi,\pi)=J\cup\{n+k\}$.
\end{proof}

Now we can prove Theorem~\ref{thm:j2}.

\begin{proof}[Proof of Theorem~\ref{thm:j2}]
  Assume that $J'=J\cup\{n+m\}$ is a $j_2$-set. Then the three largest integers
  of $J'$ are $n-k,n$ and $n+m$. If $m=1$, we are done. If $m\geq2$, then by
  Lemma~\ref{thm:2} we have $m\geq k$.

  Now assume that $m=1$ or $m\geq k$. If $m=1$, then clearly $J\cup\{n+m\}$ is a
  $j_2$-set. Assume that $m\geq k$. If $k=1$, then by Lemma~\ref{thm:5}
  $J\cup\{n+m\}$ is a $j_2$-set. If $k\geq2$, then by Lemmas~\ref{thm:3} and
  \ref{thm:4} we have that $J\cup\{n+i\}$ is a $j_2$-set for all integers
  $i\geq k$. In particular $J\cup\{n+m\}$ is a $j_2$-set. 
\end{proof}

\section{Proof of Lemma~\ref{thm:tlim}}\label{sec:pf}

In this section we prove Lemma~\ref{thm:tlim}. First we prove the following three
lemmas.

\begin{lem}\label{thm:qinv}
We have
$$\frac{t_n(q^{-1})}{\xfactorial{n}{q^{-1}}}
=\frac{t_n(q)}{\qfactorial{n}},$$
$$\frac{A_n(p^{-1},q)}{\xfactorial{n}{p^{-1}}\xfactorial{n}{q}}=
\frac{A_n(p,q^{-1})}{\xfactorial{n}{p}\xfactorial{n}{q^{-1}}}=
\frac{A_n(p^{-1},q^{-1})}{\xfactorial{n}{p^{-1}}\xfactorial{n}{q^{-1}}}=
\frac{A_n(p,q)}{\pfactorial{n}\qfactorial{n}}.$$
\end{lem}
\begin{proof}
For $T\in\t_n$ we have $\maj(T)+\maj(T')=\binom{n}{2}$,
where $T'$ is the transpose of $T$.
Thus
$$t_n(q^{-1}) = \sum_{T\in\t_n} q^{-\maj(T)} =
q^{-\binom{n}{2}}\sum_{T\in\t_n} q^{\maj(T')} =q^{-\binom{n}{2}}
t_n(q).$$ Since $q^{-\binom{n}{2}}/\xfactorial{n}{q^{-1}}
=1/\qfactorial{n}$, we get the first identity. The rest identities can be proved
similarly.
\end{proof}

\begin{lem}\label{thm:bound}
For a real number $0<q<1$, we have
$$\log \left(\prod_{i\geq1}(1-q^i)^{-i}\right) < \left( 1+
\frac{q}{(1-q)^2} \right) \left( 1+ \log \frac{1}{1-q} \right).$$
\end{lem}
\begin{proof}
The left hand side is equal to
  \begin{align*}
& \sum_{i\geq1} i\log \left(\frac{1}{1-q^i}\right)
=\sum_{i\geq1} i \sum_{j\geq1} \frac{q^{ij}}{j}\\
&=\sum_{i,j\geq2} \frac{iq^{ij}}{j} +\sum_{i\geq1} iq^i +\sum_{j\geq1}
\frac{q^j}{j} -1\\
&< \sum_{i,j\geq2} \frac{iq^{i+j}}{j} +\sum_{i\geq1} iq^i +\sum_{j\geq1}
\frac{q^j}{j}\\
&< \left(1+\sum_{i\geq1}iq^i \right)
\left(1+\sum_{j\geq1} \frac{q^j}{j}\right)
=\left( 1+ \frac{q}{(1-q)^2} \right) 
\left( 1+ \log \frac{1}{1-q} \right).
  \end{align*}
\end{proof}

Let $\x=(x_1,x_2,\ldots)$ and $\y=(y_1,y_2,\ldots)$ be two infinite
sequences of independent variables.  Let $s_{\lambda}(\x)$ denote the
Schur function in the variables $\x$.  The following
formulas are well known, see \cite{Macdonald1995, Stanley1999}:
\begin{align}\label{eq:cauchy}
\sum_{n\geq0}\sum_{\lambda\vdash n} s_{\lambda}(\x) z^{n} &= \prod_{i\geq 1}
(1-x_iz)^{-1} \prod_{1\leq i<j} (1-x_ix_j z^2)^{-1},\\ 
\sum_{n\geq0}\sum_{\lambda\vdash n}
s_{\lambda}(\x) s_{\lambda}(\y) z^{n} &= \prod_{i,j\geq 1}(1-x_i y_j z)^{-1}.
\end{align}

It is also known, see \cite[Proposition 7.19.11]{Stanley1999}, that 
$$s_{\lambda}(1,q,q^2,\ldots) = \frac{f^{\lambda}(q)}
{(1-q)^n \qfactorial{n}}.$$ Thus we get the following:
\begin{equation}\label{thm:schur1}
\sum_{\lambda\vdash n} s_{\lambda}(1,q,q^2,\ldots) =
  \frac{t_n(q)}{(1-q)^n\qfactorial{n}},
\end{equation}
\begin{equation}\label{thm:schur2}
\sum_{\lambda\vdash n} s_{\lambda}(1,p,p^2,\ldots)
  s_{\lambda}(1,q,q^2,\ldots) = \frac{A_n(p,q)}{(1-p)^n(1-q)^n
\pfactorial n \qfactorial n}.
\end{equation}

\begin{lem}\label{thm:limitlem}
For real numbers $0<p,q<1$, we have
$$\limit{n} \sum_{\lambda\vdash n}s_{\lambda}(1,q,q^2,\ldots) =
\prod_{i\geq 1} (1-q^i)^{-1} \prod_{0\leq i<j} (1-q^{i+j})^{-1},$$
$$\limit{n} \sum_{\lambda\vdash n}
s_{\lambda}(1,p,p^2,\ldots) s_{\lambda}(1,q,q^2,\ldots) =
\prod_{\substack{i,j\geq0\\ i+j>0}} (1-p^iq^j)^{-1}.$$
\end{lem}
\begin{proof}
Let $\xi_n(q)=\sum_{\lambda\vdash n}
s_{\lambda}(1,q,q^2,\ldots)$.
By \eqref{eq:cauchy}, we have
$$\sum_{n\geq0} \xi_n(q) z^n = \prod_{i\geq 0} (1-q^iz)^{-1}
\prod_{0\leq i<j} (1-q^{i+j}z^2)^{-1},$$ equivalently,
$$\sum_{n\geq0} (\xi_n(q) -\xi_{n-1}(q)) z^n = \prod_{i\geq 1}
(1-q^iz)^{-1} \prod_{0\leq i<j} (1-q^{i+j}z^2)^{-1},$$ where
$\xi_{-1}(q)=0$.

Then
$$\limit{N} \xi_N(q) = \limit{N}\sum_{n=0}^N 
(\xi_n(q) -\xi_{n-1}(q)) = \prod_{i\geq 1}
(1-q^i)^{-1} \prod_{0\leq i<j} (1-q^{i+j})^{-1}$$
converges, because
$$\prod_{i\geq1} (1-q^i)^{-1}< \prod_{0\leq i<j} (1-q^{i+j})^{-1} =
\prod_{i\geq1}(1-q^i)^{-\left\lceil \frac i2 \right\rceil} <
\prod_{i\geq1}(1-q^i)^{-i},$$ where $\prod_{i\geq1}(1-q^i)^{-i}$
converges by Lemma~\ref{thm:bound}. Thus we get the first limit.
Similarly, we can prove the second limit.
\end{proof}

\begin{proof}[Proof of Lemma~\ref{thm:tlim}]
We will only prove the first limit. The second can be proved
similarly.  Using the well known asymptotic behavior of $t_n \sim
\frac{1}{\sqrt2}n^{n/2} \exp\left(-\frac{n}{2}+\sqrt n
-\frac{1}{4}\right)$, we can easily see that it holds for
$q=1$. 

Assume $q\ne 1$.  By Lemma~\ref{thm:qinv}, it is sufficient to
show that for $0<q<1$,
$$\limit{n} \frac{\displaystyle\frac{t_{n+1}(q)}{\qfactorial{n+1}}}
{\displaystyle\frac{t_n(q)}{\qfactorial{n}}} = 1-q.$$ 
Since
$$\frac{t_n(q)}{\qfactorial{n}} = (1-q)^n \sum_{\lambda\vdash
  n}s_{\lambda}(1,q,q^2,\ldots),$$
 we are done by
Lemma~\ref{thm:limitlem}.
\end{proof}

\section{Further study}\label{sec:further-study}

In this paper we have found a $q$-analog and a $(p,q)$-analog of the following:
\begin{equation}
  \label{eq:fs1}
\limit{n} \frac{|\t_n(A)|}{|\t_n|}=
\frac{f^{\alpha}}{a!},
\end{equation}
\begin{equation}
  \label{eq:fs2}
\limit{n} \frac{|\t_n(A,B)|}{|\t_n(\emptyset,\emptyset)|} =
    \frac{f^{\alpha}f^{\beta}}{a!b!}.  
\end{equation}

In a probabilistic point of view \eqref{eq:fs1} means that the probability that
a random SYT $T$ of size $n$ contains a given SYT $A$ of shape $\alpha\vdash a$
tends to $\frac{f^\alpha}{a!}$ as $n$ approaches infinity. This clearly implies
that the probability that a random pair $(P,Q)$ of SYTs of size $n$ contains a
given pair $(A,B)$ of SYT of shape $\alpha\vdash a$ and $\beta\vdash b$ tends to
$f^\alpha f^\beta/a! b!$ as $n$ approaches infinity. However,
\eqref{eq:fs2} implies that adding the additional condition that we must have
$\sh(P)=\sh(Q)$ does not change the probability, which is nontrivial.  With this
observation we conjecture the following.

\begin{conj}\label{conj1}
  For given SYTs $A_1,\ldots,A_k$ with $\sh(A_i)=\alpha_i\vdash a_i$ for $i\in[k]$, we
  denote by $\t_n(A_1,\ldots,A_k)$ the set of $k$-tuples $(T_1,\ldots,T_k)$ such
  that $\sh(T_1)=\cdots=\sh(T_k)$ and $T_i\in \t_n(A_i)$ for all $i\in[k]$. Then
\[
\limit{n} \frac{|\t_n(A_1,\ldots,A_k)|}{|\t_n(E_1,\ldots,E_k)|} =
    \frac{f^{\alpha_1} \cdots f^{\alpha_k}}{a_1!\cdots a_k!},\]
where $E_1=\cdots=E_k=\emptyset$.
\end{conj}

\begin{problem}
  Find a $q$-analog of Conjecture~\ref{conj1} if it is true. 
\end{problem}

In our $q$-analog and $(p,q)$-analog of \eqref{eq:fs1} and \eqref{eq:fs2} we
have the assumption $p,q>0$. We can extend the ranges of $p$ and $q$ by proving
Lemma~\ref{thm:tlim} for more general $p$ and $q$. Thus we propose the following
problem.

\begin{problem}
  Find negative real numbers or complex numbers $p$ and $q$ satisfying the following:
$$\limit{n} \frac{\displaystyle \frac{t_{n+1}(q)}{\qfactorial{n+1}}}
{\displaystyle \frac{t_n(q)}{\qfactorial{n}}} = 1-\qq,\quad \limit{n}
\frac{\displaystyle \frac{A_{n+1}(p,q)}{\pfactorial{n+1}\qfactorial{n+1}}}
{\displaystyle \frac{A_n(p,q)}{\pfactorial{n}\qfactorial{n}}}
 = (1-\pp)(1-\qq).$$
\end{problem}

\section*{Acknowledgements}
I would like to thank Professor Richard Stanley for helpful discussion.  I would
also like to thank Professor Ron King for pointing out an error in
Lemma~\ref{thm:limitlem} in an earlier version of this paper. I am grateful to
the anonymous referees for their very careful reading and helpful comments which
have improved the readability of this paper significantly.

\bibliographystyle{plain}

\end{document}